\documentclass[runningheads]{lncse}
\font\elevenbb=msbm10 at 10.95pt

\def\DR{\hbox{\elevenbb DR}}
\def\K{\hbox{\elevenbb K}}
\def\N{\hbox{\elevenbb N}}
\def\M{\hbox{\elevenbb M}}
\def\R{\hbox{\elevenbb R}}

\def\Gr{Gr\"obner }
\def \bg #1 {\begin{tabular}{{#1}}}
\def \nd {\end{tabular}}
\newcommand \mwhile {{\bf while}\hspace{0.3cm}}
\newcommand \mrepeat {{\bf repeat}\hspace{0.3cm}}
\newcommand \muntil {{\bf until}\hspace{0.3cm}}
\newcommand \mfore {{\bf for\hspace{0.3cm}each}\hspace{0.3cm}}
\newcommand \mdo {{\bf do}\hspace{0.3cm}}

\newcommand \mif {{\bf if}\hspace{0.3cm}}

\newcommand \mthen {{\bf then}\hspace{0.3cm}}
\newcommand \melse {{\bf else}\hspace{0.3cm}}

\newcommand \mchoose {{\bf choose}\hspace{0.3cm}}

\newcommand \mbegin {{\bf begin}}
\newcommand \mend {{\bf end}}

\newcommand \minput {{\bf Input:\bb}}
\newcommand \moutput {{\bf Output:\bb}}
\newcommand \bb {\hspace{0.3cm}}
\newcommand \h {\hspace{0.5cm}}
\newcommand \hh {\hspace{1.0cm}}
\newcommand \hhh {\hspace{1.5cm}}
\newcommand \hhhh {\hspace{2.0cm}}
\newcommand \hhhhh {\hspace{2.5cm}}

\newcounter{cc}

\newcommand \hln {\hfill \addtocounter{cc}{1} \arabic{cc}
   \vskip -2pt \noindent }
\newenvironment{algorithm}[1]
 {
  \vskip 0.3cm
  \begin{tabular}[htb]{|p{11cm}|}
   \multicolumn{1}{c}{
  Algorithm: {\bf #1} } \setcounter{cc}{0} 
  \\  \hline
 }
 {\\ \hline \end{tabular}
  \vskip 0.3cm
  \noindent
  }
\begin{document}
\title{On the Relation Between Pommaret and Janet
Bases\thanks{Dedicated to the memory of Alyosha Zharkov.}
}
\titlerunning{On the Relation Between Pommaret and Janet
Bases}
\author{Vladimir P.Gerdt}
\authorrunning{Vladimir P.Gerdt}
\institute{Laboratory of Computing Techniques and Automation,
Joint Institute for Nuclear Research, 141980 Dubna, Russia}
\maketitle
\begin{abstract}
In this paper the relation between Pommaret and Janet bases
of polynomial ideals is studied. It is proved that if an ideal
has a finite Pommaret basis then the latter is a minimal
Janet basis. An improved version of the related algorithm for computation
of Janet bases, initially designed by Zharkov, is described. For an
ideal with a finite Pommaret basis, the algorithm computes this basis.
Otherwise, the algorithm computes a Janet basis which need not be minimal.
The obtained results are generalized to linear differential ideals.
\end{abstract}
\section{Introduction}

Pommaret and Janet bases may be cited as typical representatives
of involutive bases~\cite{GB1,GB2,Apel98,G99a}. Involutive bases
are \Gr bases, though, generally redundant, and involutive methods
provide an alternative approach to computation of \Gr bases. In so
doing, polynomial Pommaret bases which were first introduced
in~\cite{ZB93} have become a research subject in commutative
algebra. They can be considered as generalized left \Gr bases in
the commutative ring with respect to non-commutative
grading~\cite{Apel95}. Pommaret bases of homogeneous ideals in
generic position coincide with their reduced \Gr
bases~\cite{Mall}.  The use of these bases makes more
accessible the structural information of zero-dimensional
ideals~\cite{Zh96}. Pommaret bases provide an algorithmic tool for
determining combinatorial decompositions of polynomial modules
and for computations in the syzygy modules \cite{Seiler}.

Linear differential Pommaret bases form the main tool in formal theory
of linear partial differential equations~\cite{Pommaret78,Pommaret94}
whereas linear differential Janet bases form an algorithmic tool in
Lie symmetry analysis of nonlinear differential
equations~\cite{Schwarz1,Schwarz2}.
Unlike reduced \Gr bases, Pommaret and Janet bases along with any other
involutive bases lead to explicit formulae for Hilbert functions
and Hilbert polynomials~\cite{Apel98,G99a,Pommaret94,Janet}.

Basic properties of Pommaret and Janet bases are determined by the
underlying involutive divisions~\cite{GB1,GB2,Apel98}. Non-noetherity
of Pommaret division~\cite{GB1} is responsible for non-existence of
finite Pommaret bases for some polynomial (linear differential)
ideals of positive (differential) dimension. On the other hand, any
polynomial ideal as well as any linear differential ideal
has a finite Janet basis due to the noetherity of Janet division.
The two divisions differ greatly by their definition: unlike Janet
divisibility, Pommaret divisibility do not depend on the leading terms of
generators. Given an ideal and an admissible  monomial ordering, or ranking
in the differential case, its monic Pommaret basis is unique whereas there
are infinitely many different monic Janet bases and among them
only the minimal Janet basis is uniquely defined~\cite{GB2}.

However, in spite of the above differences, Pommaret and Janet
bases are closely related, and Zharkov was the first to observe
this fact in the last paper~\cite{Zh94} of his life. He argued
that if a polynomial ideal has a finite Pommaret basis and a Janet
basis which is Pommaret autoreduced, then they have identical
monic forms (c.f. \cite{Seiler}). Zharkov put also forward an
algorithm for construction of Janet bases by sequential treatment
of Janet nonmultiplicative prolongations followed by Pommaret
autoreduction.

The goal of this paper is to study the relation between polynomial
Pommaret and Janet bases in more details and to improve the Zharkov algorithm.
Our analysis is based on the properties of Janet and Pommaret divisions
and involutive algorithms studied in~\cite{GB1,GB2,G99a,G99b}.

This paper is organized as follows. The next section sketches some definitions,
notations and conventions which are used in the sequel.
Section 3 deals with analysis of the relationships between polynomial
Pommaret and Janet bases. In particular, we prove that if an ideal has a
finite Pommaret basis, then it is a minimal Janet basis.
We describe here an algorithm of the combined Pommaret and Janet autoreduction
which, given a Janet basis, converts it into another Janet basis. Since
a minimal Janet basis is both Janet and Pommaret autoreduced, the existence
of a finite Pommaret basis is equivalent to Pommaret-Janet autoreducibility
of any non-minimal Janet basis into a minimal one.
In Section 4 we describe an algorithm for computation of polynomial Janet
bases which is an improved version of the Zharkov algorithm~\cite{Zh94}. One
of the improvements is the use of Pommaret-Janet autoreduction rather then
the pure Pommaret autoreduction. Another improvement is incorporation of the
involutive analogue~\cite{GB1} of Buchberger's chain criterion~\cite{Buch85}.
Section 5 contains generalization of the results of Sections 3 and 4 to
linear differential ideals. The generalization is based on
paper~\cite{G99a} where general involutive methods and algorithms of
papers~\cite{GB1,GB2,G99b} are extended from commutative to
differential algebra.

\section{Basic Definitions and Notations}

Let $\N$ be the set of nonnegative integers, and $\M=\{x_1^{d_1}\cdots
x_n^{d_n}\ |\ d_i\in \N\}$ be a set of monomials in the polynomial ring
$\R=K[x_1,\ldots,x_n]$ over a field $K$ of characteristic zero.

By $deg(u)$ and $deg_i(u)$ we denote the total degree of $u\in \M$
and the degree of variable $x_i$ in $u$, respectively. If monomial $u$
divides monomial $v$ we shall write $u|v$. Throughout the paper we
restrict ourselves with admissible monomial orderings~\cite{BWK93}
$\succ$ which are compatible with
\begin{equation}
x_1\succ x_2\succ\cdots\succ x_n\,. \label{var_order}
\end{equation}
The leading monomial of the polynomial $f\in
\R$ with respect to $\succ$ will be denoted by $lm(f)$.
If $F\subset \R$ is a polynomial set, then by $lm(F)$ we
denote the leading monomial set for $F$, and $Id(F)$ will denote the ideal
in $R$ generated by $F$. The initial ideal of an ideal $I\in \R$ with respect to
the monomial ordering $\succ$ will be denoted by $in_{\succ}(I)$. The support
of a polynomial $f$, that is, the set of monomials occurring in $f$ with
nonzero coefficients will be denoted by $supp(f)$. For the least common multiple of two
monomials $u,v\in \M$ we shall use the conventional notation $lcm(u,v)$.


\begin{definition}\cite{Pommaret78,Janet} For a monomial
$u=x_1^{d_1}\cdots x_k^{d_k}$ with $d_k>0$ the variables $x_j,j\geq k$ are
considered as {\em Pommaret multiplicative} or {\em $P-$multiplicative} and
the other variables as {\em $P-$nonmultiplicative}. For $u=1$ all the
variables are
$P-$multiplicative.
\label{def_div_P}
\end{definition}

\begin{definition}\cite{Janet} Let $U\subset \M$ be a finite monomial
set.
For each $1\leq i\leq n$ divide $U$ into groups labeled by non-negative
integers
$d_1,\ldots,d_i$:
$$[d_1,\ldots,d_i]=\{\ u\ \in U\ |\ d_j=deg_j(u),\ 1\leq j\leq i\ \}.$$
A variable $x_i$ is called {\em Janet multiplicative} or
{\em $J-$multiplicative} for $u\in U$ if $i=1$ and
$deg_1(u)=\max\{deg_1(v)\ |\ v\in U\}$, or if
$i>1$, $u\in [d_1,\ldots,d_{i-1}]$ and $deg_i(u)=\max\{deg_i(v)\ |\ v\in
[d_1,\ldots,d_{i-1}]\}$.  \label{div_J}
If a variable is not $J-$mul\-ti\-pli\-ca\-tive
 for $u\in U$, it is called
{\em $J-$nonmultiplicative} for $u$.
\label{def_div_J}
\end{definition}

\noindent
For a polynomial $f\in \R$ the Pommaret separation of variables into
multiplicative and nonmultiplicative is done in accordance with
Definition~\ref{def_div_P} where $u=lm(f)$. Analogously, for an element
$f\in F$ in a finite polynomial set $F\subset \R$ the Janet multiplicative
and nonmultiplicative variables are determined by
Definition~\ref{div_J} with $u=lm(f)\in U=lm(F)$.

We denote by $M_P(f)$, $NM_P(f)$ and by $M_J(f,F)$, $NM_J(f,F)$,
respectively, the sets of $P-$multiplicative, $P-$nonmultiplicative
and $J-$mul\-ti\-pli\-ca\-tive, $J-$nonmultiplicative variables for $f$.
A set of monomials in $P-$mul\-ti\-pli\-ca\-tive variables for $u$ and
$J-$multiplicative variables for $u\in U$ will be denoted by $P(u)$ and
$J(u,U)$, respectively.

\begin{remark}
The monomial sets $P(u)$ and $J(u,U)$ for any $u,U$ such as $u\in U$
satisfy the following axioms
\begin{tabbing}
~~(a)~~\=If $w\in L(u,U)$ and $v|w$, then $v\in L(u,U)$. \\
~~(b)  \>If $u,v\in U$ and $uL(u,U)\cap vL(v,U)\not=\emptyset$, then \\
~~~~~  \> $u\in vL(v,U)$ or $v\in uL(u,U)$. \\
~~(c)  \> If $v\in U$ and $v\in uL(u,U)$, then $L(v,U)\subseteq L(u,U)$. \\
~~(d)  \> If $V\subseteq U$, then $L(u,U)\subseteq L(u,V)$ for all $u\in V$.
\end{tabbing}
if one takes either $P(u)$ or $J(u,U)$ as $L(u,U)$. The axioms characterize
an involutive monomial division, a concept invented in~\cite{GB1}.
Every monomial set $L(u,U)$ satisfying the axioms generates an appropriate
separation of variables into $(L-)$multiplicative and $(L-)$nonmultiplicative.
As this takes place, an element $u\in U$ is an $L-$divisor of a monomial
$w\in \M$ if $w/u\in L(u,U)$. In this case $w$ is $L-$multiple of $u$.
Using the axioms, a number of new divisions was constructed~\cite{GB2,G99b}
which may be also used for algorithmic computation of involutive bases.
\label{rem_id}
\end{remark}

\noindent
All the next definitions in this section are those in~\cite{GB1,GB2}
specified to Pommaret and Janet divisions.

\begin{definition}Given a finite monomial set $U$, its {\em cone} $C(U)$,
 {\em $P-$cone} $C_P(U)$ and {\em $J-$cone} $C_J(U)$ are the following
monomial sets
 $$ C(U)=\cup_{u\in U}\,u\M, \quad C_P(U)=\cup_{u\in U}\,uP(u),
 \quad C_J(U)=\cup_{u\in U}\,uJ(u,U)\,.$$
\label{dev_cone}
\end{definition}

\begin{definition}Given an admissible ordering $\succ$, a polynomial set
 $F\subset \R$ is called {\em Pommaret autoreduced} or {\em $P-$autoreduced}
 if every $f\in F$ has no terms $P-$multiple of an element in $lm(F)\setminus
lm(f)$.
 Similarly, a finite polynomial set $F$ is {\em Janet autoreduced} or
 {\em $J-$autoreduced} if each term in every $f\in F$ has no $J-$divisors
 among $lm(F)\setminus lm(f)$. A finite set $F$ will be called {\em
 Pommaret-Janet autoreduced} or {\em $PJ-$autoreduced} if it is both
 $P-$autoreduced and $J-$autoreduced.
\label{aut_PJ}
\end{definition}

\begin{remark}
From Definition~\ref{div_J} it follows that any finite set $U$ of distinct
monomials is $J-$autoreduced
$$
(\forall u,v\in U)\ (u\neq v)\ \ [\ uJ(u,U)\cap vJ(v,U)=\emptyset\ ]\,.
$$
\label{J_aut_ms}
\end{remark}

\begin{definition}Given an admissible ordering $\succ$ and a polynomial
 set $F\subset \R$, a polynomial $h\in \R$ is said to be in {\em $P-$normal
form
 modulo $F$} if every term in $h$ has no $P-$ divisors in $lm(F)$.
 Similarly, if all the terms in $h$ have no  $J-$divisors among the leading
terms
 of a finite polynomial set $F$, then $h$ is in {\em $J-$normal form modulo
$F$}.
\label{def_nf}
\end{definition}

\noindent
A general involutive normal form algorithm is described in~\cite{GB1},
and an involutive normal form of a polynomial modulo any involutively
autoreduced
set is uniquely defined. We denote by $NF_P(f,F)$ and $NF_J(f,F)$,
respectively,
$P-$normal and $J-$normal form of polynomial $f$ modulo $F$.

Pommaret or Janet autoreduction of a finite polynomial set $F$ may
be performed~\cite{GB1} similar to the conventional
autoreduction~\cite{Buch85,BWK93}. If $H$ is obtained from $F$ by
the conventional, or $J-$autoreduction we shall write $H=Autoreduce(F)$
or $H=Autoreduce_J(F)$, respectively.

\begin{definition}A $P-autoreduced$ set $F$ is called
 a {\em Pommaret basis} ($P-$basis) of the ideal $Id(F)$ generated by $F$ if
\begin{equation}
 (\forall f\in F)\ (\forall x\in NM_P(f))\ \ [\ NF_P(f\cdot x,F)=0\
]\,.\label{Pb}
\end{equation}
Similarly, a $J-$autoreduced set $F$ is called a {\em Janet basis}
($J-$basis)
 of $Id(F)$ if
\begin{equation}
 (\forall f\in F)\ (\forall x\in NM_J(f,F))\ \ [\ NF_J(f\cdot x,F)=0\ ]\,.
\label{Jb}
\end{equation}
\label{def_PJ_bases}
\end{definition}

\noindent
In accordance with Definition \ref{def_nf} the nonmultiplicative prolongation $f\cdot x$
with the vanishing $P-$ or $J-$normal form modulo polynomial set
$F=\{f_1,\cdots,f_m\}$ can be rewritten as
$$
f\cdot x=\sum_{i=1}^m f_i\,h_i,\quad h_i\in \R,\quad lm(f)\cdot x\succeq
lm(f_i\,h_i)\quad (1\leq i\leq m)
$$
where $supp(h_i)\subset P(f_i)$ or $supp(h_i)\subset J(f_i,F)$,
respectively, for every polynomial product $f_i\,h_i$.

Let $G_P$ and $G_J$ be Pommaret and Janet bases of an ideal
$I$, respectively. Then, from Definition~\ref{def_PJ_bases} it
follows~\cite{GB1} that
\begin{equation}
h\in I\quad \mbox{iff}\quad NF_P(h,G_P)=0\,\quad \mbox{and}\quad
h\in I\quad \mbox{iff}\quad NF_J(h,G_J)=0\,. \label{id_mem}
\end{equation}
This implies, in particular, the equalities
\begin{equation}
C_P(lm(G_P))=C(lm(G_P))\,,\quad C_J(lm(G_J))=C(lm(G_J))\,.\label{con_ims}
\end{equation}

\noindent
It is immediate from (\ref{con_ims}) that
$lm(G_P)$ and $lm(G_J)$ are $P-$and $J-$bases of the initial ideal
$in_{\succ}(I)$.

\begin{corollary} If for a $P-$autoreduced set $G_P$ the equality
(\ref{con_ims}) of its cone and $P-$cone holds and $lm(G_P)$ is a
basis of the initial ideal $in\left(Id(G_P)\right)\}$, then
$G_P$ is a $P-$ basis of $Id(G_P)$. Analogous statement holds for
a $J-$basis. \label{cr_lm}
\end{corollary}

\noindent
Whereas monic Pommaret bases much like to reduced \Gr bases
are unique, every ideal, by property (\ref{J_aut_ms}), has
infinitely many monic Janet bases. Among them there is the unique
$J-$basis defined as follows.

\begin{definition}A Janet basis $G$ of ideal $Id(G)$ is called {\em
minimal}
if for any other $J-$basis $F$ of the ideal the inclusion $lm(G)\subseteq
lm(F)$ holds.
\label{def_minJb}
\end{definition}

\begin{remark}Every zero-dimensional polynomial ideal has a finite
Pommaret
basis, and for a positive dimensional ideal the existence of finite Pommaret
basis
can be always achieved by means of an appropriate linear transformation of
variables~\cite{ZB93,Apel95,Pommaret78,Pommaret94}.
\label{finite_Pb}
\end{remark}

\section{Relation Between Polynomial Pommaret and Janet Bases}

Given a finite monomial set $U$, Definitions~\ref{def_div_P} and
\ref{def_div_J} generally give different separations of variables for
elements in $U$.

\begin{example} $U=\{x_1x_2,x_2x_3,x_3^2\}$.
\begin{center}
\vskip 0.2cm
\bg {|c|c|c|c|c|} \hline\hline
monomial &\multicolumn{2}{c|}{Pommaret}
& \multicolumn{2}{c|}{Janet} \\ \cline{2-5}
         &   $M_P$     & $NM_P$    & $M_J$         & $NM_J$  \\ \hline
$x_1x_2$ &   $x_2,x_3$ & $x_1$     & $x_1,x_2,x_3$ & $-$   \\
$x_2x_3$ &   $x_3$     & $x_1,x_2$ & $x_2,x_3$     & $x_1$    \\
$x_3^2$  &   $x_3$     & $x_1,x_2$ & $x_3$         & $x_1,x_2$  \\ \hline
\hline
\nd
\vskip 0.2cm
\end{center} \label{exm_1}
\end{example}

\noindent
Here is, however, an important relation between Pommaret and Janet
separations:

\begin{proposition}{\em\cite{Zh94}(see also \cite{GB1})}. If a finite
monomial set
$U$ is $P-$auto\-re\-duced, then for any $u\in U$ the following inclusions hold
$$M_P(u,U)\subseteq M_J(u,U),\quad NM_J(u,U)\subseteq NM_P(u,U)\,.$$
\label{id_P_J}
\end{proposition}

\noindent For $U$ in Example~\ref{exm_1} the minimal Janet basis $U_J$ and
the Pommaret basis $U_P$ of the monomial ideal $Id(U)$ are
\begin{equation}
U_J=U\cup \{x_1x_3^2\}, \quad U_P=U_J\cup_{i=2}^\infty
\{x_1^ix_2\}\cup_{j=2}^\infty \{x_1^jx_3^2\} \cup_{k=2}^\infty \{x_2^kx_3\}\,.
\label{inf_PB}
\end{equation}
Note, that $U_P$ is
infinite and $U_J\subset U_P$. Below we show that the inclusion
$G_J\subseteq G_P$ holds for any minimal Janet basis $G_J$ and Pommaret basis
$G_P$ if both of them are monic and generate the same ideal.
Furthermore, the proper inclusion $G_J\subset G_P$ holds iff $P$ is
infinite.

The following algorithm, given a finite polynomial set $F\in \R$
and an admissible ordering $\succ$, performs $PJ-$auto\-re\-duc\-tion of
$F$ and outputs a $PJ$-\-auto\-re\-duced set $H$. In this case we shall
write $H=Autoreduce_{PJ}F$.

Since the involutive $P-$ and $J-$reductions which are performed in the
course of the algorithm form subsets of the conventional
reductions~\cite{GB1}, the algorithm terminates. Furthermore, the
{\bf while}-loop generates the $P-$auto\-re\-duced monomial set
$lm(H)$. In accordance with Remark \ref{J_aut_ms} the Janet
autoreduction in line 12 does not affect the leading terms,
and, hence,  the output polynomial set is both Pommaret and Janet
autoreduced.

\begin{center}
\begin{algorithm}{Pommaret-JanetAutoreduction($F,\,\succ$)} \label{PJA}
\minput $F\in \R$, a finite set; $\succ$, an admissible ordering \\
\moutput $H=Autoreduce_{PJ}(F)$ \\
\mbegin \\
\h $G:=\emptyset$;\ \ $H:=F$
\hln
\h \mrepeat
\hln
\hh $\tilde{H}:=H$
\hln
\hh \mwhile exist $h\in H$ such that
\hln
\hhh $lm(h)\in C_P\left(lm(H\setminus \{h\})\right)$\bb \mdo
\hln
\hhh $H:=H\setminus \{h\}$;\ \ $G:=G\cup \{h\}$
\hln
\hh \mfore $g\in G$\bb \mdo
\hln
\hhh $G\setminus \{g\}$;\ \ $f:=NF_J(g,H)$
\hln
\hhh \mif $f\neq 0$\bb \mthen
\hln
\hhhh $H:=H\cup \{f\}$
\hln
\h \muntil $H=\tilde{H}$
\hln
\h $H:=Autoreduce_J(H)$
\hln \mend
\end{algorithm}
\end{center}

\begin{proposition}
If algorithm {\bf Pommaret-JanetAutoreduction} takes a Janet basis $F$ as an input
its output $H$ is also a Janet basis of the same ideal, and $H\subseteq F$.
\label{prop_PJ_aut}
\end{proposition}

\begin{proof}
Let $F$ be a Janet basis of the ideal $Id(F)$ and
$H$ be a polynomial set which computed by the algorithm.
Apparently, $Id(H)=Id(F)$. Consider $U=lm(F)$. If $U$ is
$P-$autoreduced, then $H$ initiated as $F$ in line 1 does not change
in the process of the algorithm, and, hence, $H=F=Autoreduce_{PJ}(F)$.
Otherwise, consider the output polynomial sets $H$.
Denote $lm(H)$ by $V$. Then, by construction, $V\subset U$.

Consider a monomial $t\in C_J(U)$ and show that $t\in C_J(V)$.
Let $u\in U$ be a $J-$divisor of $t$, that is, $t\in uJ(u,U)$, and
$v\in V$ be such that $v|u$. If $u=v$, by property
(d) of Janet division in Remark \ref{rem_id}, we are done. Let
now $u\in U\setminus V$. The {\bf while}-loop provides that $u\in vP(v)$.
If $v=x_1^{d_1}\cdots x_k^{d_k}$ with $d_k>0$ $(1\leq k\leq n)$,
then from Definition \ref{def_div_P} it follows that
$u=x_1^{d_1}\cdots x_k^{d_k+e_k}\cdots x_n^{e_n}$.

We have to prove that any variable $x_i\in J(u,U)$
$(1\leq i\leq n)$ satisfy $x_i\in J(v,V)$. Consider two
alternative cases: $i<k$ and $i\geq k$. In the first case,
by Definition \ref{def_div_J}, both $u,v$ belong to the
same group $[d_1,\cdots,d_{i-1}]$ of monomials in $U$.
It follows that $x_i\in J(v,V)$. In the second case, by
Definition \ref{def_div_P}, $x_i\in P(v)$ and Proposition
\ref{id_P_J} we find again that $x_i$ is $J-$multiplicative
for $v$ as an element in $V$.

Therefore, $V$ is a Janet monomial basis of $Id(U)$. Thus, by Corollary \ref{cr_lm},
$H$ is a $J-$basis.  In so doing, every $J-$normal form computed in
line 8 of the algorithm apparently vanishes. This implies the inclusion
$H\subseteq F$.
\qed
\end{proof}

\begin{corollary}
A minimal Janet basis is Pommaret autoreduced.
\label{cr_mJb}
\end{corollary}

\begin{proof}
Let $G$ be a minimal Janet basis. From
Proposition~\ref{prop_PJ_aut} and Definition~\ref{def_minJb}
it follows that $lm(G)$ is $P-$autoreduced. Thus, $G$ is $P-$autoreduced.
\qed
\end{proof}

\noindent
It is clear that, given a Janet basis, its $PJ-$autoreduction
yields, generally, more compact basis than the pure
$P-$autoreduction.

\begin{example} Consider polynomial set
$$ F=\{xyzt-xz,xyz+z^2,xzt+x^2,xy+z,zt+x\}$$
which is a Janet basis of the ideal $Id(F)$ with respect to the
degree-reverse-lexicographical ordering $\succ$ such that
$x\succ y\succ z\succ t$. Given $F$ and $\succ$ as input, the
algorithm {\bf Pommaret-JanetAutoreduction}$(F,\succ)$ outputs the
minimal $J-$basis
$$G_1=\{xzt+x^2,xy+z,zt+x\}\,.$$
If one uses the Pommaret normal form computation in line 7 instead of
that of Janet, then the algorithm leads to the output
$$G_2=\{xzt+x^2,z^2t+xz,xy+z,zt+x\}=G_1\cup \{z^2t+xz\}\,.$$
\end{example}

\noindent
The following theorem is the main theoretical result of the present
paper and forms a basis of an algorithm for construction of
Janet and Pommaret bases described in the next section.

\begin{theorem} Given an ideal $I\subseteq R$ and an admissible
monomial ordering $\succ$ compatible with (\ref{var_order}), the
following are equivalent:
\begin{tabbing}
~~$(i)$~~~\=$I$ has a finite Pommaret basis. \\
~~$(ii)$~~\>A minimal Janet basis of $I$ is its Pommaret basis. \\
~~$(iii)$~\>If $F$ is a Janet basis of $I$, then $G=Autoreduce_{JP}(F)$ is
\\
~~~~~~~~\>its Pommaret basis.
\end{tabbing}
\label{th_JP_bases}
\end{theorem}

\begin{proof}
$(i)\Longrightarrow(iii)$: Suppose
$G=\{g_1,\ldots,g_m\}$ which, by Proposition \ref{prop_PJ_aut},
is also a $J-$basis of $I$, is not its $P-$basis. Our goal is to prove
that a Pommaret basis of $I$ is an infinite polynomial set.
From Corollary~\ref{cr_lm} it follows that
\begin{equation}
(\exists g\in G)\ (\exists x_k\in NM_P(g))\ \ [\ lm(g)\cdot x_k \not\in
C_P\left(lm(G)\right)\ ]\,. \label{assumption}
\end{equation}
Among nonmultiplicative prolongations $g\cdot x_k$
satisfying (\ref{assumption}) choose one with the lowest $lm(g)\cdot x_k$
with respect to the pure lexicographical ordering $\succ_{Lex}$ generated
by (\ref{var_order}). If there are
several such prolongations choose that with the lexicographically lowest
$x_k$, that is, with the lexicographically highest $g$.
This choice is unique since $G$ is $J-$autoreduced.

We claim that $x_k\in M_J(g,G)$. Assume for a contradiction that $x_k\in
NM_J(g,G)$. Then from Janet involutivity conditions (\ref{Jb}) we obtain
\begin{equation}
(\exists f\in G)\ \ [\ lm(g)\cdot x_k=lm(f)\cdot w\
|\ w\in J\left(lm(f),lm(G)\right)\ ]\,. \label{ass_cond}
\end{equation}
In accordance with condition (\ref{assumption}) $w$ contains $P-$nonmultiplicative
variables for $f$ and from (\ref{ass_cond}) it
follows \cite{GB1} that $f\succ_{Lex}g$. If $w=x_j$, then
$x_j\succ x_k$, and both $lm(f)\cdot x_j$ and $lm(g)$ belong to
the same monomial group $[\ deg_1(g),\ldots,deg_{k-1}(g)\ ]$ appearing in
Definition \ref{def_div_J}. Hence, $x_j\in NM_J(f,G)$
in contradiction to (\ref{ass_cond}). Therefore, $deg(w)\geq 2$ and
(\ref{ass_cond}) can be rewritten as
$$ lm(g)\cdot x_k=\left(lm(f)\cdot v\right)\cdot x_j\, $$
with $w=v\cdot x_j$.
Show that $lm(f)\cdot v\in P\left(lm(f)\right)$. If
$v\in C_P\left(lm(G)\right)$ we are done. Otherwise there is $x_m|v$
such that $x_m\in NM_P\left(lm(f)\right)$. Since
$lm(f)\cdot x_m\prec_{Lex} lm(g)\cdot x_k$, our choice of $g$ and $x_k$
implies
the existence $g_1\in G$ such that $f\cdot x_m = g_1\cdot v_1$
and $v/x_m \in P\left(lm(g_1)\right)$. If $v_1=v/x_m\in P\left(lm(g_1)\right)$
we are done. Otherwise we select again a $P-$nonmultiplicative
variable for $g_1$ occurring in $v_1$ and rewrite the corresponding prolongation
in terms of its $P-$divisor $lm(g_2)\in lm(G)$. Continuity of Pommaret
division \cite{GB1} provides termination of the rewriting process with
an element in $\tilde{g}\in G\setminus \{g\}$ such that $lm(\tilde{g})$ is a $P-$divisor
of $lm(f)\cdot v$. Because $G$ is $P-$autoreduced, by Proposition
\ref{id_P_J} $lm(\tilde{g})$ is also a $J-$divisor of $lm(f)\cdot v$.
By this means there are two different Janet divisors $lm(f)$
and $lm(\tilde{g})$ of the same monomial that contradicts Remark \ref{J_aut_ms}
and proves the claim.

Let now $H$ be a $P-$basis of $Id(G)$. Denote $lm(g)\cdot x_k$ by $u$
and show that $u\in lm(H)$. Suppose there is an element $h\in H\setminus G$,
such that $u=lm(h)\cdot v$ with $v\in P\left(lm(h)\right)$. Then
$h\prec_{Lex} u$ and there is $q\in G$ satisfying $lm(h)=lm(q)\cdot w$
where $w\not \in P\left(lm(q) \right)$.
Thus, there is a $P-$nonmultiplicative prolongation $q\cdot x_j$ with
$x_j|w$ such that $lm(q)\cdot x_j\prec_{Lex} u$ and
$lm(q)\cdot x_j\not \in C_P\left(lm(G)\right)$ that contradicts
the above choice of $g$ and $x_k$.

Now consider monomial set $U=lm(G)\cup \{u\}\subseteq lm(H)$. By Definition
\ref {def_div_J}, $x_k\in NM_P(u)$ and $u\cdot x_k$ is obviously the
lexicographically lowest $P-$nonmultiplicative prolongation of elements in $U$.
It is easy to see that $u\cdot x_k\not \in C_P(U)$. Indeed, since $x_k\in M_J(g,G)$,
it follows that $u\cdot x_k\not \in U$, and a $P-$divisor of $u\cdot x_k$ would
also $P-$divide $lm(g)$ that is impossible as $G$ is $P-$autoreduced. Thus, we
find that $u\cdot x_k\in lm(H)$. By sequential repetition of this
reasoning for $u\cdot x_k^i$ $(i=2,3\ldots)$ we deduce that every such
monomial is an element in $lm(H)$, and, therefore, $H$ is infinite.

$(iii)\Longrightarrow (ii)$: If $G$ is a minimal Janet basis of $Id(G)$,
then Corollary \ref{cr_mJb} implies $G=Autoreduce_{JP}(G)$, and the above
arguments show that either $G$ is also a $P-$basis or the
latter is infinite.

$(ii)\Longrightarrow (i)$: This implication is obvious.
\qed
\end{proof}

\begin{corollary} Let $G_J$ be a monic minimal Janet basis and $G_P$ be a monic
Pommaret basis for the same polynomial ideal and monomial
ordering. Then $G_J\subseteq G_P$ and $G_J\subset G_P$ iff $G_P$
is infinite.
\label{cr_subbasis}
\end{corollary}

\begin{proof} This follows from the above proof of Theorem
\ref{th_JP_bases}.
\qed
\end{proof}


\begin{corollary} If a polynomial ideal is in generic position, then
its minimal Janet basis is also a Pommaret basis. If such an ideal is
homogeneous, then these bases are also reduced \Gr bases of the ideal.
\label{cr_generic}
\end{corollary}

\begin{proof} $I$ has a finite Pommaret basis \cite{Apel95,Mall}. If
$I$ is homogeneous, then its Pommaret basis is the reduced \Gr basis
\cite{Mall}.
\qed
\end{proof}

\section{Algorithm for Construction of Janet Bases}

In this section we present an algorithm for constructing $J-$bases
of polynomial ideals which is based on Theorem~\ref{th_JP_bases}
and will be called {\bf JanetBasis}. Whenever the
ideal generated
by an input polynomial set has a finite $P-$basis for a given
admissible ordering, the algorithm outputs just
this basis which is also a minimal $J-$basis. Otherwise,
the $J-$basis computed by the algorithm is not necessarily minimal
as we demonstrate below by the explicit example.

\begin{center}
\begin{algorithm}{JanetBasis($F,\,\succ$)} \label{JB}
\minput $F\in \R$, a finite set; $\succ$, an admissible ordering \\
\moutput $G$, an involutive basis of the ideal $Id(F)$ \\
\mbegin \\
\h $G:=Autoreduce(F)$
\hln
\h $T:=\emptyset$
\hln
\h \mfore $g\in G$\bb \mdo
\hln
\hh $T:=T\cup \{(g,lm(g),\emptyset)\}$
\hln
\h \mwhile exist $(g,u,D)\in T$ and $x\in NM_J(g,G)\setminus D$
  \bb \mdo
\hln
\hh \mchoose such $(g,u,D),x$ with the lowest
 $lm(g)\cdot x$ w.r.t. $\succ$
\hln
\hh $T:=T\setminus \{(g,u,D)\} \cup \{(g,u,D\cup \{x\})\}$
\hln
\hh \mif $Criterion(g\cdot x,u,T)$ is false \bb \mthen
\hln
\hhh $h:=NF_J(g\cdot x,G)$
\hln
\hhh \mif $h\neq 0$\bb \mthen
\hln
\hhhh \mif $lm(h)=lm(g\cdot x)$\bb \mthen
\hln
\hhhhh $T:=T\cup \{(h,u,\emptyset)\}$
\hln
\hhhh \melse
\hln
\hhhhh $T:=T\cup \{(h,lm(h),\emptyset)\}$
\hln
\hhhhh $G:=Autoreduce_{PJ}(G\cup \{h\})$
\hln
\hh $Q:=T$
\hln
\hh $T:=\emptyset$
\hln
\hh \mfore $g\in G$\bb \mdo
\hln
\hhh \mif exist $(f,u,D)\in Q$ s.t. $lm(f)=lm(g)$\bb \mthen
\hln
\hhhh \mchoose $\tilde{g}\in G$ s.t.
$u\in lm(\tilde{g})J\left(lm(\tilde{g}),lm(G)\right)$
\hln
\hhhh $T:=T\cup \{\left(g,lm(\tilde{g}),D\cap NM_J(g,G)\right)\}$
\hln
\hhh \melse
\hln
\hhhh $T:=T\cup \{(g,lm(g),\emptyset)\}$
\hln
\mend
\end{algorithm}
\end{center}
\noindent
$Criterion(g,u,T)$ is true provided that if there is $(f,v,D)\in T$
such that $lcm(u,v) \prec lm(g)$ and $lm(g)\in lm(f)J\left(lm(f),lm(G)\right)$.

The structure of algorithm {\bf JanetBasis} is very close to that of
the algorithm {\bf InvolutiveBasis} devised in~\cite{GB1} for
construction of involutive bases for arbitrary constructive
involutive divisions, and, hence, applicable to Janet division.
The main difference between the algorithms is the form of their
intermediate autoreduction.
Whereas the previous algorithm, when specified for Janet division,
uses the pure Janet autoreduction which do not affect the leading terms,
the below one uses the above algorithm {\bf Pommaret-JanetAutoreduction}.
By this reason, the algorithm {\bf InvolutiveBasis}, unlike the below one,
almost never outputs a minimal Janet basis or a Pommaret basis if the
latter is finite. In paper~\cite{GB2} we designed the algorithm
{\bf MinimalInvolutiveBasis}
which always outputs a minimal involutive basis whenever the latter exists.
As we now see from Theorem \ref{th_JP_bases}, this algorithm
in the case of Janet division outputs also a Pommaret basis
if it is finite.

Besides, in the algorithm {\bf InvolutiveBasis} the involutive
autoreduction is performed whenever nonzero normal form is obtained
in the {\bf while}-loop. In the algorithm {\bf JanetBasis} the
subalgorithm {\bf Pommaret-JanetAuto\-re\-duc\-tion}
is caused by a nonzero $J-$norm\-al form $h$ in line 15 only if the
leading
term of the related prolongation is $J-$reducible. Otherwise, since
$G$ is always $P-$autoreduced before its enlargement with $h$, $lm(h)$ cannot
$P-$divide, in accordance with Definition \ref{def_div_P}, any other
element in $lm(G)$.

Note that we indicated the intersection in line 21 to emphasize that
elements in $D$ must be nonmultiplicative variables for the corresponding
polynomial that is always understood in algorithms of papers
(\cite{GB1,GB2,G99a,G99b}).

Noetherity of Janet division provides termination of the algorithm
{\bf JanetBasis} \cite{GB1}. To show this
consider the intermediate bases $G_0=Autoreduce(F)$ and $G_i$ $(i=1,2,\ldots)$
generated after the $i-th$ iteration of the {\bf while}-loop. It is
clear that
\begin{equation}
Id\left(lm(G_0)\right)\subseteq Id\left(lm(G_1)\right)\subseteq
Id\left(lm(G_2)\right)\subseteq \cdots \label{chain}
\end{equation}
and this chain is stabilized after finitely many steps. Namely, the
stabilization starts when the intermediate polynomial set
$G$ becomes a (non-necessarily reduced) \Gr basis
of $Id(F)$. By partial involutivity of $G$ \cite{GB1}, the proper
inclusion in chain
(\ref{chain}) holds only when $G$ is enlarged by $h$ in line 15. In between
of such proper inclusions and after the chain stabilization
$lm(G)$ is completed with $lm(h)=lm(g\cdot x)$ as
stands in line 11, and this completion cannot be infinite by noetherity
of Janet division \cite{GB1}. Once algorithm terminates, it produces,
by Proposition \ref{prop_PJ_aut},
a $PJ-autoreduced$ Janet basis of $F$ because the involutivity conditions
(\ref{Jb})
are satisfied as is checked in lines 6 and 9 where $h=0$. In so doing,
correctness of the criterion which is verified in line 8 is proved
exactly as it done in \cite{GB1} for the algorithm {\bf InvolutiveBasis}.

\begin{remark} In line 6 of the algorithm {\bf JanetBasis} one
can use any admissible ordering for selection of the current
$J-$nonmultiplicative prolongation $g\cdot x$ to be treated in the
following lines. This selection ordering may not only be different from
the main ordering $\succ$ but also may vary at every step when the
selection is done. Correctness of this arbitrariness in the choice of
selection ordering follows from the correctness of this arbitrariness for
the monomial completion procedure \cite{G99b}.
 \label{rm_order}
\end{remark}

\noindent
As mentioned above, the algorithm {\bf JanetBasis} for an ideal of positive
dimension may not output its minimal $J-$basis. We demonstrate this by the
following example.

\begin{example}Let $F$ be a set $\{x^2y-z,x\,y^2-y\}$ generating
one-dimensional ideal, and $\succ$ be
the degree-reverse-lexicographical ordering with $x\succ y\succ z$. Then
the algorithm {\bf JanetBasis$(F,\succ$)} outputs the following Janet
basis
$$ G=\{x^2y-z,x^2z-z^3,x\,y-y\,z,x\,z-z^2,y^2z-y,y\,z^2-z\} $$
whereas the minimal Janet basis coinciding with the reduced \Gr basis is
$$ \{x\,y-y\,z,x\,z-z^2,y^2z-y,y\,z^2-z\}\,. $$
\label{ex_2}
\end{example}

\begin{remark} The algorithm {\bf JanetBasis} is an improved version of
the algorithm designed by Zharkov~\cite{Zh94}. The first improvement is
the use of the mixed Pommaret-Janet autoreduction instead of the pure
Pommaret autoreduction as Zharkov proposed. Let $G_1$ and $G_2$ be output
Janet bases computed with the use of $PJ-$ and $P-$autoreduction,
respectively. Then, Proposition \ref{id_P_J} implies $G_1\subseteq G_2$
and below we give an example when the proper inclusion holds. The second
improvement is the criterion used in line 7. This criterion is an
involutive analogue~\cite{GB1} of the Buchberger's chain criterion~\cite{Buch85}
and is superior to the criterion used in~\cite{Zh94}.
\label{rem_zharkov}
\end{remark}

\begin{example}\cite{Karin} The following polynomial set generates
three-dimensional ideal
\begin{equation}
F=\left\{
\begin{array}{l}
4x_5 (x_6 + 2a_1-8x_1)(a_2-a_3)-x_2x_3x_4+x_2+x_4 \\
4x_5 (x_6 + 2a_1-8x_2)(a_2-a_3)-x_1x_3x_4+x_1+x_3 \\
4x_5 (x_6 + 2a_1-8x_3)(a_2-a_3)-x_1x_2x_4+x_2+x_4 \\
4x_5 (x_6 + 2a_1-8x_4)(a_2-a_3)-x_1x_2x_3+x_1+x_3
\end{array}
\right\}  \label{speer}
\end{equation}
where
$$a_1=x_1+x_2+x_3+x_4,\ a_2=x_1x_2x_3x_4,\ a_3=x_1x_2+x_2x_3+x_3x_4+x_4x_1\,.$$
For the polynomial set (\ref{speer}) and the degree-reverse-lexicographical
ordering compatible with (\ref{var_order}) the algorithm {\bf Janet basis}
outputs set $G_1$ with 71 polynomials. The pure Pommaret autoreduction in line
15 generates set $G_2\supset G_1$ of 75 elements. Note that a minimal Janet
basis contains 49 polynomials whereas the reduced \Gr basis contains 44
polynomials.
\label{ex_speer}
\end{example}

\section{Linear Differential Bases}

Let $\K$ be a zero characteristic differential field
with a finite number of mutually commuting derivation operators
$\partial /\partial x_1,\ldots,\partial /\partial x_n$. Consider
the differential polynomial ring $\DR=\K\{y_1,\ldots,y_m\}$ with the
set of differential indeterminates $\{y_1,\ldots,y_m\}$.
Elements in $\DR$ are differential polynomials in
$\{y_1,\ldots,y_m\}$. An ideal in $\DR$ generated by linear
differential polynomials is called linear differential
ideal~\cite{Kolchin}.

In~\cite{G99a}, by exploiting the well-known algorithmic similarities
between polynomial and linear differential systems and the association
between monomials and the derivatives, we extended the general involutive
methods and algorithms designed in \cite{GB1,GB2,G99b} for polynomial ideals
to linear differential ideals. The statements and algorithms of Sect. 3
and 4 admit similar extension. As this takes place, all the
above statements have proven for the polynomial case are proved
by parallel arguments for the differential case. In the following
table we give a short correspondence between these two cases. In
particular, this correspondence allows one to rewrite the algorithms
{\bf Pommaret-JanetAutoreduction} and {\bf JanetBasis} for linear
differential bases.

\begin{center}
\vskip 0.2cm
\bg {|c|c|} \hline\hline
 & \\[-0.2cm]
\ Commutative Algebra \  &\ Differential Algebra \ \\
 & \\[-0.2cm] \hline
 &  \\
$\ \R=K[x_1,\ldots,x_n]\ $   & $\ \DR=\K\{y_1,\ldots,y_m\}\ $    \\
$ f,g,h \in \R $ & $ f,g,h \in \DR $ \\
$ F,G,H \subset \R $ & $ F,G,H\subset \DR $ \\
$ x^\alpha\equiv x_1^{\alpha_1}\cdot x_n^{\alpha_n}$ &\ $\partial_{\alpha}
y_j \equiv \frac {\partial^{\alpha_1+\cdots+\alpha_n} y_j}
{\partial x_1^{\alpha_1} \cdots \partial x_n^{\alpha_n}} \Longleftrightarrow [x^\alpha]_j\ $ \\
$g\cdot x$ & $\partial_x g$ \\
monomial ordering  & ranking \\
$lm(g)$  & leading derivative $ld(g) $ \\
$lm(G)$  & $ld(G)=\cup_{g\in G} \{ld(g)\}$ \\
 &
\\ \hline
\hline
\nd
\vskip 0.2cm
\end{center}

\section{Conclusion}

As we have seen finite Pommaret bases of polynomial and linear
differential ideals are minimal Janet bases. The above proof of Theorem
\ref{th_JP_bases} shows that, given a $P-$autoreduced Janet basis $G$,
it is easy to verify the existence of a finite Pommaret basis. One
suffices to check the condition (\ref{assumption}). Another check
which may be even easier in practice is to look at the Pommaret and
Janet separation of variables for elements in $G$. As shown in
\cite{Seiler}, the existence of a finite Pommaret basis implies
coincidence of both separations for every element in $G$. Moreover,
the leading monomial structure of an infinite Pommaret basis can be
read off the structure of $lm(G)$ (c.f. (\ref{inf_PB}) for Example \ref{exm_1}).

Therefore, minimal Janet bases can be used in commutative and differential
algebra as well as a Pommaret bases with the advantage of finiteness. For
example, in the formal theory of differential
equations~\cite{Seiler,Pommaret78,Pommaret94} infinity of a Pommaret basis
signals on $\delta$-singularity of the coordinate system chosen, and the
condition (\ref{assumption}) gives the same signal for Janet bases.

On the face of it, Pommaret division looks like more attractive than
Janet one since its separation is easier to compute than the Janet separation.
Besides, Pommaret division, unlike that of Janet, is globally defined \cite{GB2}.
Hence, after enlargement of the intermediate polynomial set with an
irreducible $P-$nonmultiplicative prolongation there are no needs to recompute
the Pommaret separation for other elements in the set. However, the careful
implementation of both divisions do not reveal any notable advantage of Pommaret
division over Janet division in construction of polynomial bases. This rather
surprising fact was firstly observed by Zharkov \cite{Zh94}. One of the explanations of this
experimental phenomenon is given by Proposition \ref{id_P_J}: one must generally
treat more $P-$nonmultiplicative prolongations than $J-$nonmultiplicative ones.
In so doing, the search for an involutive divisor in the process of involutive
reduction can be done similarly for both divisions as we show in our
forthcoming paper \cite{GBY}.

The algorithm {\bf JanetBasis} presented above is now under implementation in
C in parallel with the algorithm {\bf MinimalInvolutiveBasis} \cite{GB2}
specified for Janet division. Our first experimentation with the codes
shows that sometimes the former algorithm runs faster than the latter one and
needs less computer memory. For example, modular computation of the
degree-reverse-lexicographic Janet basis for the Cyclic 7 example~\cite{BM} is
about twice faster with the algorithm {\bf JanetBasis} than with the
algorithm {\bf MinimalInvolutiveBasis}. Currently, the timings for computation
modulo 31013 are about 5 and 10 minutes, respectively, on a Pentium Pro 333
Mhz computer.

\section{Acknowledgements} I am grateful to W.M.Seiler for useful
discussions and for providing me with his recent paper \cite{Seiler}.
I would also like to thank D.A.Yanovich for his assistance in computing
the above examples. This work was partially supported by grant INTAS-99-1222
and grants No. 98-01-00101 and 00-15-96691 from the Russian Foundation for
Basic Research.

\end{document}